\documentclass[a4paper,UKenglish]{lipics}
 
\usepackage{microtype}

\title{The Roles Leon Henkin Played in Mathematics Education}
\titlerunning{Henkin on Math Education} 

\author{Mar\'{\i}a Manzano}
\affil{Salamanca University, Department of Philosophy\\
  Campus Unamuno, Edificio FES, 37007 Salamanca\\
  \texttt{mara@usal.es}}

\authorrunning{M. Manzano} 

\Copyright{María Manzano}


\serieslogo{logo_ttl}
\volumeinfo
  {M. Antonia {Huertas}, Jo\~ao {Marcos}, Mar\'ia {Manzano}, Sophie {Pinchinat}, \\
  Fran\c{c}ois {Schwarzentruber}}
  {5}
  {4th International Conference on Tools for Teaching Logic}
  {1}
  {1}
  {111}
\EventShortName{TTL2015}

\begin{document}

\maketitle

\begin{abstract}
This paper is divided into two sections. In the first I give reasons for
strongly recommending reading some of Henkin's expository papers. In the
second I describe  Leon Henkin's work as a social activists in the field of
mathematics education, as he labored in much of his career to boost the
number of women and underrepresented minorities in the upper echelons of
mathematics.
\end{abstract}

\section{Some of Henkin's expository papers}

Henkin was an extraordinary insightful professor with a talent for
exposition, and he devoted considerable effort to writing expository papers.
I will mention three of them, while trying to convince you to read them with
your students as a source of mutual inspiration.

\subsection{Are Logic and Mathematics Identical?}

This is the title of a wonderful expository paper \cite{1962a}, which Leon
Henkin published in \emph{Science} in 1962, and subtitled: \emph{An old
thesis of Russell's is reexamined in the light of subsequent developments in
mathematical logic.}

I recommend that you to give this paper to your students, not only because
the historical view provided is comprehensive and synthetic but also because
it shows the Henkin's characteristic style; namely, the ability to strongly
catch your attention from the start.

\emph{How does he achieve it? }you might wonder. In that particular paper,
Henkin tells us that his interest in logic began at the age of 16, when 
\emph{`(I)\ came across a little volume of Bertrand Russell entitled \textbf{%
Mysticism and Logic}'. } In the introduction Henkin cites Russell's radical
thesis, that\emph{\ `mathematics was nothing but logic'} together with the
companion thesis \emph{\ ` that logic is purely tautological'}, and he
describes the strong reaction against his thesis by the academic community: 
\emph{`Aux armes, citoyens du monde math\'{e}matique!'}

Henkin then devotes the first section of the paper to explain the two main
ideas that could help explain how Russell arrived at his conclusion. The
first  was the lengthy effort to achieve a \emph{`systematic reduction of
all concepts of mathematics to a small number of them'}, and the second was 
\emph{`the systematic study by mathematical means of the laws of logic which
entered into mathematical proofs'. }Henkin relates the work of Frege,
Peirce, Boole and Schr\"{o}der, during the second half of the nineteenth
century, with the two efforts mentioned above, and identifies them as the
primary \emph{raison d'etre }of \emph{Principia Mathematica.}

In the following section, entitled \emph{From Russell to G\"{o}del}, Henkin
explains the introduction of semantic notions by Tarski, as well as the
formulation and proof of the completeness theorem for propositional logic by
Post and for first-order logic by G\"{o}del. \emph{`This result of G\"{o}%
del's is among the most basic and useful theorems we have in the whole
subject of mathematical logic'. }But, Henkin also explains how, in 1931, the
hope of further extension of this kind of completeness was \emph{`dashed by G%
\"{o}del himself [...] (he) was able to demonstrate that the system of 
\textbf{Principia Mathematica}, taken as a whole was incomplete'.\ }%
Immediately after, and anticipating what the reader might be thinking,
Henkin dispels the hope of finding new axioms to repair the incompleteness
phenomenon.

In the section entitled \emph{Consistency and the Decision Problem}, Henkin
analyzes these important notions and also explains how \emph{`G\"{o}del was
able to show that the questions of consistency and completeness were very
closely linked to one another. [...] if a system such as the \textbf{%
Principia} were truly consistent, then in fact it would not be possible to
produce a sound proof of this fact!'. }In the following section, named \emph{%
Logic after 1936}, Henkin describes how Alonzo Church proved that no
decision procedure is available for first-order logic, and he devotes the
rest of the paper  to set theory, recursive functions, and algebraic logic.
Henkin ends the paper with a section where he analyzes \emph{Russell's
Thesis in Perspective}.\emph{\ }

\bigskip

Henkin was awarded the \emph{Chauvenet Prize} in 1964 for this paper. The
prize is described as a \emph{Mathematical Association of America} award to
the author of an outstanding expository article on a mathematical topic by a
member of the Association.

\subsubsection{Bertrand Russell's request}

In April 1, 1963, Henkin received a very interesting letter from Bertrand
Russell. In it, Russell thanked Henkin for \emph{`your letter of March 26
and for the very interesting paper which you enclosed.' }Right at the
beginning Russell declared:

\begin{quote}
It is fifty years since I worked seriously at mathematical logic and almost
the only work that I have read since that date is G\"{o}del's. I realized,
of course, that G\"{o}del's work is of fundamental importance, but I was
pussled by it. It made me glad that I was no longer working at mathematical
logic. If a given set of axioms leads to a contradiction, it is clear that
at least one of the axioms is false. Does this apply to school-boys'\
arithmetic, and if so, can we believe anything that we were taught in youth?
Are we to think that 2+2 is not 4, but 4.001?
\end{quote}

He then went on explaining his \emph{`state of mind'} while Whitehead and he
were doing the Principia and added: \emph{`Both Whitehead and I were
disappointed that the Principia was almost wholly considered in connection
with the question whether mathematics is logic.'}

Russell ended the letter with a request: \emph{`If you can spare the time, I
should like to know, roughly, how, in your opinion, ordinary mathematics
---or, indeed, any deductive system--- is affected by G\"{o}del's work.'}

According to Annellis: \emph{`Henkin replied to Russell at length with an
explanation of G\"{o}del's incompleteness results, in a letter of July 1963,
specifically explaining that G\"{o}del's showed, not the inconsistency, but
the incompleteness of the [Principia] system.'}

\subsection{On Mathematical Induction}

In a personal communication Henkin affirmed that \emph{On mathematical
induction }\cite{Henkin60}, published in 1960, was the favorite among his
articles because it had a somewhat panoramic nature and was not directed
exclusively to specialists. He wrote: \emph{`}[...]\emph{\ but my little
paper on induction models from 1960, which has always been my favorite among
my expository papers'}. In it, the relationship between the induction axiom
and recursive definitions is studied in depth.

\emph{Why do I so strongly recommend that you ask your students to read this
paper? }From my point of view, it is the best paper on logic to offer
students as a first reading of a \textquotedblleft
real-life\textquotedblright\      article. The paper is especially interesting
because Henkin describes something that would never appear in a formal
article: his motivation. It seems that Henkin was trying to convince a
mathematical colleague about why a given argument about the existence of
recursive operations was completely wrong, even though at first sight it
might seem convincing.

Before going into the details of the wrong argument that motivated the whole
paper, Henkin explains Peano arithmetic.

Peano axiomatized the theory of natural numbers. To do so, he started out
from indefinable primitive terms ---in particular, those of natural number,
zero, and the successor function---- and by means of three axioms he
synthesized the main facts. Among those axioms is that of induction, which
states that any subset of natural numbers, closed by the successor operation
and to which zero belongs, is precisely the set of all natural numbers.
Although the axioms for the theory of natural numbers are very important,
the most interesting theorems of the theory did not stem from them alone
because in most of the theorems, operations of addition, multiplication,\
etc. are used.

Peano thought that after axiomatizing a theory it did not suffice to
organize the facts by means of axiomatic laws; it was also necessary to
organize the concepts, using definition laws. In particular, we can define
addition recursively by means of:

\begin{enumerate}
\item $x+0=x$

\item $x+Sy=S(x+y)$, for all $x,y$ of $%
\mathbb{N}
$
\end{enumerate}

However, this definition must be justified by a theorem in which the
existence of a unique operation that will satisfy the previous equations
must be established.

A poor argument to prove this is as follows. Let us choose any $x\in 
\mathbb{N}
$. We define a subset $G$ of $%
\mathbb{N}
$, by placing all $y\in 
\mathbb{N}
$ elements in $G$ for which $x+y$ is defined by the previous equations. It
is not difficult to see that $0\in G$ and that $\forall (y\in G\rightarrow
Sy\in G)$. Now using the induction axiom, we conclude that $G=%
\mathbb{N}
$ and thus that $x+y$ is defined for all $y\in 
\mathbb{N}
$.

\emph{Why is this proof wrong? }This was the question that Henkin's
colleague posed him. Henkin tried to convince him that because the argument
was designed to establish the existence of a function $f$ ($+$ in the
example), it is incorrect to assume in the course of the argument that we
have such a function. Henkin made him see that in the proof only the third
axiom was used and that, if correct, the same reasoning could be used not
only for models that satisfy all Peano axioms but also for those that
satisfy only the induction one. Henkin called these \textquotedblleft
Induction Models\textquotedblright\ and proved that in them not all
recursive operations are definable. For example, exponentiation fails.

Induction models turn out to have a fairly simple mathematical structure:
there are standard ones ---that is, isomorphic to natural numbers--- but
also non-standard ones. The latter also have a simple structure: either they
are cycles, in particular $\mathbb{Z}$ modulo $n$, or they are what Henkin
calls \textquotedblleft spoons\textquotedblright\; because they have
a handle followed by a cycle. The reason is that the induction axiom is
never fulfilled alone, since it requires Peano's first or second axiom. This
does not mean that Peano's axioms are redundant, as it is well known that
they are formally independent; \emph{i.e.}, each one is independent of the
other two.

\subsection{Completeness}

If you take a look at the list of documents Leon Henkin left us, the first
published paper, \emph{The completeness of first order logic }\cite%
{TiposHenkin:49}, corresponds to his well known result, while the last, 
\emph{The discovery of my completeness proof }\cite{henkin1996},\emph{\ } is
a extremely interesting as autobiography, thus ending his career with a sort
of fascinating loop.

\emph{I claim that reading the last paper is a must. Why?} As you know, Leon Henkin left us an important collection of papers, some of
them so exciting as his proof of the completeness theorem both for the
theory of types and for first-order logic. He did so by means of an
innovative and highly versatile method, which was later to be used in many
other logics, even in those known as non-classical. In his 1996 paper, we
learn about the process of discovery, which observed facts he was trying to
explain, and why he ended up discovering things that were not originally the
target of his enquiries. Thus, in this case we do not have to engage in
risky hypotheses or explain his ideas on the mere basis of the later, cold
elaboration in scientific articles. It is well known that the \emph{logic of
discovery} differs from what is adopted on organizing the final exposition
of our research through their different propositions, lemmas, theorems and
corollaries.

We also learn that the publication order of his completeness results (\cite%
{TiposHenkin:49} and \cite{TiposHenkin:50}) is the reverse of his discovery
of the proofs. The completeness for first-order logic was accomplished when
he realized he could modify the proof obtained for type theory in an
appropriate way. We consider this to be of great significance, because the
effort of abstraction needed for the first proof (that of type theory)
provided a broad perspective that allowed him to see beyond some prejudices
and to make the decisive changes needed to reach his second proof. In \cite%
{manzanoCompHenkin2014} you can find a detailed commentary of Henkin's
contribution to the resolution and understanding of the completeness
phenomena.

\subsubsection{Henkin's expository papers on completeness}

In 1967 Henkin published two very relevant expository papers on for the subject we are
considering here, \emph{Truth and Provability} and \emph{Completeness},
which were published in \emph{Philosophy of Science Today} \cite%
{Morgenbesser1967}.

\paragraph{\emph{Truth and Provability}}

In less than 10 pages, Henkin gives a very intuitive introduction to the
concept of truth and its counterpart, that of provability, in the same
spirit of Tarski's expository paper \emph{Truth and Proof }\cite{Tarski1969}. The latter was published in \emph{Scientific American} two years
after Henkin's contribution. This not so surprising as Henkin had by then
been in Berkeley working with Tarski for about 15 years and the theory of
truth was Tarski's contribution.

The main topics Henkin introduces (or at least touches upon) are very
relevant. They include the \emph{use/mention} distinction, the desire for 
\emph{languages with infinite sentences} and the need for a \emph{recursive
definition of truth}, the \emph{language/metalanguage} distinction, the need
to avoid reflexive paradoxes, the concept of \emph{denotation }for terms,
and the interpretation of \emph{quantified formulas}. He also explains what
an \emph{axiomatic} \emph{theory} is and how it works in harmony with a 
\emph{deductive calculus}. Properties such as \emph{decidability} and \emph{%
completeness/incompleteness of a theory} are mentioned at the end. I admire
the way these concepts are introduced, with such \'{e}lan, and the chain
Henkin establishes, which shows how each concept is needed to support the
next.

\paragraph{\emph{Completeness}}

In this short expository paper Henkin explores the complex landscape of the
notions of completeness. He introduces the notion of logical completeness
---both weak and strong--- as an extension of the notion already introduced
of \textquotedblleft completeness of an axiomatic theory\textquotedblright .
This presentation differs notably from the standard way these notions are
introduced today where, usually, the completeness of the logic precedes the
notion of completeness of a theory and, often, to avoid misunderstandings, 
both concepts are separated as much as possible, as if relating them were
some sort of terrible mistake or even anathema. G\"{o}del's incompleteness
theorem is presented, as well as its negative impact on the search for a
complete calculus for higher-order logic. The paper ends by introducting his
own completeness result for higher-order logic with general semantics. The
utilitarian way Henkin uses to justify his general models as a way of
sorting the provable sentences from the unprovable ones in the class of
valid sentences (in standard models) is very peculiar.

\section{The Roles of Action and Thought in Mathematics Education}

Henkin was often described as a social activist, he labored much of his
career to boost the number of woman and underrepresented minorities in the
upper echelons of mathematics. He was also very aware that we are beings
immersed in the crucible of history from which we find it hard to escape, an
awareness he brought to the very  beginning of his interesting article
about the teaching of mathematics \cite{henkin1995}:

\begin{quote}
Waves of history wash over our nation, stirring up our society and our
institutions. Soon we see changes in the way that all of us do things,
including our mathematics and our teaching. These changes form themselves
into rivulets and streams that merge at various angles with those arising in
parts of our society quite different from education, mathematics, or
science. Rivers are formed, contributing powerful currents that will produce
future waves of history.

The Great Depression and World War II formed the background of my years of
study; the Cold War and the Civil Rights Movement were the backdrop against
which I began my career as a research mathematicians, and later began to
involve myself with mathematics education.\newline
(In \cite{henkin1995} page 3)
\end{quote}

In this paper he gave both a short outline of the variety of educational
programs he created and/or participated in, and interesting details about
some of them. In particular he discussed the following six:

\begin{enumerate}
\item 1957-59, \textbf{NSF} \textsc{Summer Institutes}. The National Science
Foundation is an independent federal agency created by Congress in 1950. As
you can read in their web page, \emph{http://www.nsf.gov/about/}, its aim
was \emph{`to promote the progress of science; to advance the national
health, prosperity, and welfare; to secure the national defense\ldots '}.
Nowadays, NSF is the \emph{`only federal agency whose mission includes
support for all fields of fundamental science and engineering, except for
medical sciences.'} NSF's Strategic Plan includes Investing in Science,
Engineering, and Education for the Nation's Future. In \cite{henkin1995}
Henkin related this initiative to historical facts: \emph{`The launching of
Sputnik demonstrated superiority in space travel, and our country responded
in a variety of ways to improve capacity for scientific and technical
developments'} In 1957, Henkin involved himself in several NSF's programs,
he served as a lecturer of several courses. These programs were designed to
improve high school and college mathematics instruction and were directed to
mathematics teachers. Henkin explained that the variety of attitudes toward
mathematics of the teachers attending the courses was amazing, and that the
experience gave him a view of the nature of instruction around the country.
The subject of his courses was the \emph{axiomatic foundation of number
systems}. One of his aims was to get students to understand
\textquotedblleft the idea of a proof\textquotedblright\ because he believed
that it could help students in the effort of finding proofs of their own, in
a much better way than the mere understanding of the steps that constitutes
a proof.

\item 1959-63, \textbf{MAA} \textsc{Math Films.} The \textbf{M}athematical 
\textbf{A}ssociation of \textbf{A}merica was established one century ago, in
1915. As you can read in their web page, \emph{http://www.maa.org/}, \emph{%
`Over our first century, MAA has certainly grown, but continues to maintain
our leadership in all aspects of the undergraduate program in mathematics'.}
Long before internet resources became available, the MAA made movies. As
Henkin said: \emph{`Sensing a potential infusion of technology into
mathematics instruction, MAA set up a committee to make a few experimental
films. }[...]\emph{\ the committee approached me in 1959-60 with a request
to make a filmed lecture on mathematical induction which could be shown at
the high-school-senior/college-freshman level. I readily agreed.' }The film
was part of the Mathematics Today series, and was shown on public television
in New York City and in high schools. In \cite{henkin1995} Henkin explained
the preparation of the film, both from a technical point of view and from a
methodological and pedagogical perspective. He attributed the lack of
understanding of the induction principle at the undergraduate level to the
current formulation as a mathematical principle, and he proposed to use it
as \emph{`a statement about sets of numbers satisfying two simple
conditions; formulated in this way, it is a fine vehicle for giving students
practice in forming and using sets of numbers to show that all natural
numbers possess various properties'}

\item 1961-64 \textsc{CUPM}. The Mathematical Association of America's 
\textbf{C}ommittee on the \textbf{U}nder- graduate \textbf{P}rogram in 
\textbf{M}athematics (CUPM) is charged with making recommendations to guide
mathematics departments in designing curricula for their undergraduate
students. In the sixties, the CUPM proposed courses to be taken by
elementary teachers. In \cite{henkin1995} Henkin said \emph{`Some of my
colleagues and I began, for the first time, to have classroom contact with
prospective elementary teachers, and that led, in turn, to in-service
programs for current teachers. I learned a great deal from teaching
teachers-students; I hope they learned at least half as much as I!'}

\item 1964-. \textsc{Activities To Broaden Opportunity}. \textquotedblleft
The sixties\textquotedblright\ is the term used to describe the
counterculture and revolution movement that took place in several places in
the U.S.A. and Europe. Berkeley students were taking energic actions against
segregation in southeastern U.S.A. as well as against military actions in
Vietnam. In \cite{henkin1995} Henkin said \emph{`In the midst of this
turmoil I joined in forming two committees at Berkeley which enlarged the
opportunity of minority ethnic groups for studying mathematics and related
subjects. }[...]\emph{\ We noted that while there was a substantial black
population in Berkeley and the surrounding Bay Area, our own university
student body was almost \textquotedblleft lily white\textquotedblright\ and
the plan to undertake action through the Senate was initiated' }In 1964,
Leon Henkin and Jerzy Neyman, a world-famous Polish-American statistician
from Berkeley University, started a program at Berkeley to increase the
number of minority students entering college from Bay Area high schools.
Henkin told us that the inicitiave came after Neyman participation in \emph{%
`the MAA's Visiting Mathematician Program in Fall 1963. He lectured in
southerns states where, by law, whites and blacks studied in separate
colleges. Upon returning to Berkeley he told some of his friend that
\textquotedblleft first-rate students were being given a third-rate
education\textquotedblright\ ' }Henkin and Neyman undertake actions through
the Senate, and in 1964 the Senate established a committee with the desired
effect. The committee recruited promising students and offered them summer
programs to study mathematics and English. If they persisted in the program,
they were offered special scholarships. \newline
In the same year, 1964, Henkin heard a talk by a Berkeley High School
teacher, Bill Johntz. After that, Henkin was invited to see him in action,
while he was teaching mathematics to elementary students from low-income
neighborhoods, and realized that Johntz was able to raise great enthusiasm
in the class. Significantly, students enjoyed and actively engaged in the
process of learning, and they became integrally involved in their own
education. He was using a Socratic group-discovery method modeled after the
filmed teaching of David Page, a University of Illinois mathematics
professor. The method was working well, and they recruited university
mathematics students as well as engineers as teachers, after some training.
The program was called Project \emph{SEED} ---Special Elementary Education
for the Disadvantaged. This program is still alive, as you can see in their
web page, \emph{http://projectseed.org/.}

\item 1960-68. \textsc{Teaching Teachers, Teaching Kids}. In this paper
Henkin described several conferences on school mathematics as well as
several projects and courses he was involved in. The following paragraph
caught my eye: \emph{`After I began visiting elementary school classes in
connection with CTFO, I came to believe that the emotional response of the
teachers to mathematics was of more importance to the learning process of
the students than the teacher's ability to relate the algorithms of
arithmetic to the axioms of ring theory'. }

\item 1968-70. \textsc{Open Sesame: The Lawrence Hall Of Science.} The
Lawrence Hall of Science, a science museum in Berkeley, was created in honor
of the 1939 Novel prize winner Ernest Orlando Lawrence. As you can read in
their web page,%
\[
\emph{http://www.lawrencehallofscience.org/about} 
\]%
\emph{`We have been providing parents, kids, and educators with
opportunities to engage with science since 1968.'}

According to Henkin's tale \emph{`In 1968, the newly appointed director,
Professor of Physics Alan Portis, decided to transform the museum into a
center of science and mathematics education, whose functions would be
integrated with graduate research programs directed by interdisciplinary
group of faculty.'} To help in his endeavor, he gathered a group of faculty
from a variety of science departments interested in science education. \emph{%
`These faculty members proposed a new, interdisciplinary Ph.D program under
the acronym SESAME ---Special Excellence in Science and Mathematics
Education. Entering students were required to have a masters degree in
mathematics or in one of the sciences. Courses and seminars in theories of
learning, cognitive science, and experimental design were either identified
in various departments, or created'. }Nitsa Movshovitz-Hadar, a student from
the Technion in Israel, was admitted in the SESAME program, she wrote her
thesis under the direction of Leon Henkin. Nitsa is one of the contributors
of the book, \textit{The Life and Work of Leon Henkin: Essays on His
Contributions}.
\end{enumerate}

\end{document}